\documentclass{amsart}

\usepackage{amssymb}
\usepackage{graphicx}
\newtheorem{lemma}{Lemma}[section]
\newtheorem{proposition}[lemma]{Proposition}
\newtheorem{theorem}[lemma]{Theorem}
\newtheorem{corollary}[lemma]{Corollary}
\newtheorem{problem}[lemma]{Problem}

\theoremstyle{definition}
\newtheorem{definition}[lemma]{Definition}
\newtheorem{remark}[lemma]{Remark}
\newtheorem{example}[lemma]{Example}

\newcommand{\Z}{\mathbb Z} 
\newcommand{\R}{\mathbb R} 
\newcommand{\C}{\mathbb C} 
\newcommand{\Q}{\mathbb Q} 

\newcommand{\M}{\mathcal{M}} 
\newcommand{\I}{\mathcal{I}} 
\begin{document}

\title[$L^2$-torsion invariants and the Magnus representation]
{$L^2$-torsion invariants and the Magnus representation 
of the mapping class group}


\author[T. Kitano and T. Morifuji]
{Teruaki Kitano and Takayuki Morifuji}

\dedicatory{Dedicated to Professor 
Shigeyuki Morita on his 60th birthday}

\thanks{This research was partially supported by Grant-in-Aid for 
Scientific Research (No.17540064 and No.17740032), 
the Ministry of Education, Culture, Sports, Science and 
Technology, Japan.}

\address{Department of Information Systems Science, 
Soka University, 
Tokyo 192-8577, Japan}
\address{Department of Mathematics, 
Tokyo University of Agriculture and Technology, 
Tokyo 184-8588, Japan}




\begin{abstract}
In this paper, 
we study 
a series of $L^2$-torsion invariants from the viewpoint 
of the mapping class group of a surface. 
We establish some vanishing theorems for them. 
Moreover 
we explicitly calculate the first two invariants 
and compare them with hyperbolic volumes. 
\end{abstract}

\maketitle



\section{Magnus representation}

Let $\Sigma_{g,1}$\ be a compact oriented smooth surface of genus $g$ 
with a boundary $\partial \Sigma_{g,1}\cong S^1$. 
In this paper, 
we always assume that $g\geq 1$. 
We take and fix a base point
$*\in \partial \Sigma_{g,1}$ of $\Sigma_{g,1}$.
Let 
$\M_{g,1}$ be the mapping class group of
$\Sigma_{g,1}$, 
namely, 
the group of all isotopy classes of 
orientation preserving diffeomorphisms of $\Sigma_{g,1}$ 
relative to the boundary. 
We denote $\pi_1(\Sigma_{g,1},*)$ by $\Gamma$, 
which is a free group  of rank $2g$, 
and fix a generating system 
$\Gamma=\langle x_1,\dots,x_{2g}\rangle$. 
Let 
$\Z\Gamma$ be the group ring of $\Gamma$ over $\Z$. 
We write 
$\varphi_*\in \mathrm{Aut}(\Gamma)$ 
to the automorphism induced 
from $\varphi\in\M_{g,1}$. 
The following result, 
usually called the Dehn-Nielsen-Baer theorem, 
is classical and fundamental to study 
the mapping class group $\M_{g,1}$ 
by using combinatorial group theories 
(see \cite{Ivanov02-1} Section 2.9). 

\begin{proposition}[Zieschang \cite{ZVC80-1}]
The above induced homomorphism 
$\M_{g,1}$
$\ni\varphi\mapsto\varphi_*\in\mathrm{Aut}(\Gamma)$ 
is injective.
\end{proposition}

As a corollary, 
we see that 
$\varphi$ can be determined by the words 
$\varphi_*(x_1),$$\dots,$$\varphi_*(x_{2g})$$\in \Gamma$. 
Since 
the fundamental formula 
$
\gamma
=1+\sum_{i=1}^{2g}(\partial \gamma/\partial x_{i})(x_i-1)$ 
holds in $\Z\Gamma$ 
for any $\gamma\in\Gamma$, 
the word $\varphi_*(x_j)$ 
is determined by 
$\{ {\partial\varphi_{*}(x_j)}/{\partial x_i} \}$. 
Here 
$\partial/\partial x_i:\Z\Gamma\to\Z\Gamma$ 
denotes Fox's free differential. 
See 
\cite{Birman} Section 3.1 
for a systematic treatment of the subject. 
The Magnus representation of 
the mapping class group 
is defined as follows. 

\begin{definition}
The Magnus representation of $\M_{g,1}$ 
is defined by the assignment 
$$
r:\M_{g,1}\ni\varphi
 \mapsto
 \left(
 \frac{\overline{\partial\varphi_{*}(x_j)}}{\partial x_i}
 \right)_{ij}
 \in GL(2g,\Z\Gamma),
$$
where 
$\overline{\sum_g\lambda_gg}=\sum_g{\lambda_g}g^{-1}$ 
for any element 
$\sum_g\lambda_gg \in\Z\Gamma$. 
\end{definition}

\begin{remark}
By the expression 
$\gamma=1+\sum_i(\partial\gamma/\partial x_{i})(x_i-1)$, 
it follows that $r$ is injective. 
However, it is not a group homomorphism, just a crossed homomorphism. 
According to the practice, 
we call it simply the Magnus representation of $\M_{g,1}$. 
\end{remark}

Now 
for a matrix $B\in M(n,\C)$, 
let us recall that its characteristic polynomial 
$$
\mathrm{det}(tI-B)
$$ 
is one of the fundamental tools in the linear algebra. 
Here 
$I$ denotes the identity matrix of degree $n$. 
If we can define a characteristic polynomial of $r(\varphi)$, 
it may be useful tool to study the mapping class group. 
In order to 
define it for a Magnus matrix 
$r(\varphi)$, 
we need to clarify 
the following two points. 
\begin{enumerate}
\item
What is the determinant over a non-commutative group ring?
\item
What is the meaning of a variable $``t"$ in the group?
\end{enumerate}
As an answer 
to these problems, 
we can formulate that 
\begin{itemize}
\item
the variable $t$ lives in the fundamental group of the mapping torus of $\varphi$, 
\item
a characteristic polynomial 
$\mathrm{``det"}(tI-r(\varphi))$ 
with respect to 
the Fuglede-Kadison determinant. 
\end{itemize}

In the later sections, 
we explain that the characteristic polynomial 
of $r(\varphi)$ is defined as a real number 
and it essentially gives the $L^2$-torsion and 
the hyperbolic volume of the mapping torus of $\varphi$. 
Moreover 
taking the lower central series of 
the surface group $\Gamma$, 
we obtain a family of Magnus representations, 
so that 
we can introduce a sequence of $L^2$-torsion invariants 
as an approximate sequence of the hyperbolic volume. 

This paper is organized as follows. 
In the next section, 
we briefly recall the definition of the Fuglede-Kadison 
determinant. 
In Section 3, 
we summarize some properties of the $L^2$-torsion 
of 3-manifolds and 
explain a relation to the Magnus representation. 
We introduce a sequence of $L^2$-torsion invariants 
for a surface bundle over the circle 
in Section 4 and give some formulas for them 
in Section 5. 
In the last section, 
we discuss some vanishing theorems for 
$L^2$-torsion invariants. 

\section{Fuglede-Kadison determinant}

In this section, 
we review the combinatorial definition of the Fuglede-Kadison determinant 
over a non-commutative group ring  
and its basic properties 
(see \cite{Luck02-1} for details). 

The idea to define a determinant over a group ring 
comes from the following observation. 
That is, 
for a matrix $B\in GL(n,\C)$ 
with 
the (non-zero) eigenvalues 
$\lambda_1,\ldots,\lambda_n$, 
we can formally calculate 
\begin{align*}
\log |\det B|^2
&=\log\prod_{i=1}^n\lambda_i\overline{\lambda}_i
=\sum_{i=1}^n\log \lambda_i\overline{\lambda}_i\\
&=\sum_{i=1}^n\left(\sum_{p=1}^\infty\frac{(-1)^{p+1}}{p}\left(\lambda_i\overline{\lambda}_i-1\right)^p\right)\\
&=-\sum_{p=1}^\infty\left(\sum_{i=1}^n\frac1p\left(1-\lambda_i\overline{\lambda}_i\right)^p\right)\\
&=-\sum_{p=1}^\infty\frac1p\mathrm{tr}\left(\left( I-BB^*\right)^p\right)
\end{align*}
by the power series expansion of $\log$, 
where 
$B^*$ is the adjoint matrix of $B$. 
More precisely, 
if we take a sufficiently large constant $K>0$, 
we obtain 
$$
|\det B|
=K^n\exp \left(-\frac12
\sum_{p=1}^\infty\frac{1}{p}\text{tr}
\left(\left(I-K^{-2}BB^*\right)^p\right)\right)
\in {\Bbb R}_{>0}.
$$
Thus 
if we can define a certain ``trace" over a group ring, 
we get a determinant by using this formula. 

Let $\pi$ be a discrete group and 
$\C \pi$ denote its group ring over $\C$. 
For an element $\sum_{g\in\pi}\lambda_g g\in \C\pi$, 
we define the $\C\pi$-trace 
$\text{tr}_{\C\pi}:\C\pi\rightarrow \C$ 
by 
$$
\text{tr}_{\C\pi}\left(\sum_{g\in\pi}\lambda_g g\right)
=\lambda_e\in \C, 
$$
where $e$ is the unit element in $\pi$. 
For 
an $n\times n$-matrix $B=(b_{ij})\in M(n,\C\pi)$, 
we extend the definition of $\C\pi$-trace by means of 
$$
\text{tr}_{\C\pi}(B)
=\sum_{i=1}^{n}\text{tr}_{\C\pi}(b_{ii}).
$$

Next 
let us recall the definition of the $L^2$-Betti number 
of an $n\times m$-matrix $B\in M(n,m,\C\pi)$. 
We consider the bounded $\pi$-equivariant operator 
$$
R_B:\oplus_{i=1}^{n}l^2(\pi)\to\oplus_{i=1}^{m}l^2(\pi)
$$
defined by the natural right action of $B$. 
Here 
$l^2(\pi)$ is the complex Hilbert space of the formal sums 
$\sum_{g\in\pi}\lambda_g g$ which are square summable. 
We fix a positive real number $K$ so that 
$K\geq ||R_B||_\infty$ holds, 
where $||R_B||_\infty$ is the operator norm 
of $R_B$. 

\begin{definition}
{\rm 
The $L^2$-Betti number of a matrix $B\in M(n,m,\C\pi)$ is defined by 
$$
b(B)=
\lim_{p\to\infty}\ 
\text{tr}_{\C\pi}\left(\left(I-K^{-2} BB^*\right)^p\right)
\in \mathbb R_{\geq 0},
$$
where 
$B^*=(\overline{b}_{ji})$ 
and 
$\overline{\sum\lambda_g g}
=\sum\overline{\lambda}_g g^{-1}$ 
for each entry. }
\end{definition}

Roughly speaking, 
the $L^2$-Betti number $b(B)$ measures the size of the kernel of 
a matrix $B$. 
Hereafter 
we assume $b(B)=0$. 
Then, 
for a matrix with coefficients in a non-commutative group ring, 
we can introduce the desired determinant as follows. 

\begin{definition}
{\rm 
The Fuglede-Kadison determinant 
of a matrix $B\in M(n,m,\C\pi)$ is defined by 
$$
\text{det}_{\C\pi}(B)
=K^n\exp \left(-\frac12
\sum_{p=1}^\infty\frac{1}{p}\text{tr}_{\C\pi}
\left(\left(I-K^{-2} BB^*\right)^p\right)\right)
\in \mathbb R_{> 0},
$$
if the infinite sum of non-negative real numbers 
in the above exponential converges to a real number.
}
\end{definition}

\begin{remark}
It is shown that the $L^2$-Betti number $b(B)$ 
and the Fuglede-Kadison determinant $\text{det}_{\C\pi}(B)$ 
are independent of the choice of the constant $K$ 
(see \cite{Luck94-1} for example). 
\end{remark}

Here we consider the condition of the convergence. 
For any matrix $B\in M(n,\C)$, 
the condition 
$$
\lim_{p\to\infty}\ 
\text{tr}\left(\left(I-K^{-2} BB^*\right)^p\right)=0
$$ 
implies 
that $B$ has no zero eigenvalues, 
and 
then $|\det B|$ converges. 
In the case of 
group rings, 
if $\text{det}_{\C\pi}(B)$ converges, then 
$b(B)=0$. 
But 
it is not a sufficient condition, 
so that 
we need additional one. 
It is a problem to decide 
when 
$\text{det}_{\C\pi}(B)$ converges. 
Under 
the assumption that $b(B)=0$, 
such a sufficient condition is given by 
the positivity of the Novikov-Shubin invariant 
$\alpha(B)$. 
Then 
the convergence of the infinite sum 
in the Fuglede-Kadison determinant 
is guaranteed. 
The Novikov-Shubin invariant 
of an operator $R_B$ measures 
how concentrated the spectrum 
of $R_B^*R_B$ is. 
However, in general, 
it is hard to check the positivity of 
the Novikov-Shubin invariant. 

To avoid the difficulty, 
we need to consider the determinant class condition 
for groups 
(see \cite{Luck02-1}, \cite{Schick01-1} for details). 
A group $\pi$ is of $\det\geq 1$-class 
if 
for any $B\in M(n,m,\Z \pi)$ 
the Fuglede-Kadison determinant of 
$B$ satisfies $\det_{\C\pi}(B)\geq 1$. 
There are 
no known examples of groups 
which are not of $\det\geq 1$-class. 
Further 
recently it was proved that 
there is a certain large class $\mathcal{G}$ of groups 
for which 
they are of $\det\geq 1$-class. 
It includes amenable groups and 
countable residually finite groups. 
If 
we can see that $\pi$ belongs to $\mathcal{G}$, 
namely it is of $\det\geq 1$-class, 
the convergence of the Fuglede-Kadison determinant 
is guaranteed 
when the $L^2$-Betti number is vanishing. 
See \cite{Luck94-3}, \cite{Luck02-1}, \cite{Schick01-1} 
for definitions and properties 
of these subjects. 

\section{$L^2$-torsion of 3-manifolds}

In this section, 
we quickly 
recall the definition of the $L^2$-torsion of 3-manifolds. 
It is an $L^2$-analogue of the Reidemeister 
and the Ray-Singer torsion 
and essentially gives Gromov's simplicial volume 
under certain general conditions 
\cite{BFKM96-1}, 
\cite{CFM97-1}, 
\cite{CM92-1}, 
\cite{HS98-1}, 
\cite{Lott92-1}, 
\cite{LL95-1}, 
\cite{LR91-1}, 
\cite{LS99-1}, 
\cite{Mathai92-1}. 
See \cite{Luck02-1} and its references 
for historical background, 
related works and so on. 

Let $M$ be a compact connected orientable 3-manifold. 
We fix a $CW$-complex structure on $M$. 
We may assume that 
the action of $\pi_1M$ on the universal covering 
$\widetilde{M}$ is cellular 
(if necessary, we have only to take a subdivision of 
the original structure). 
We 
consider the $\C\pi_1M$-chain complex 
\[
0 \longrightarrow C_3(\widetilde{M},\C)
\overset{\partial_3}{\longrightarrow} C_2(\widetilde{M},\C)
\overset{\partial_2}{\longrightarrow} C_1(\widetilde{M},\C)
\overset{\partial_1}{\longrightarrow} C_0(\widetilde{M},\C)
\longrightarrow 0
\]
of $\widetilde{M}$. 
Since 
the boundary operator $\partial_i$ is 
a matrix with coefficients in $\C\pi_1M$, 
if we take the adjoint operator 
$\partial_i^*:C_{i-1}(\widetilde{M},\C)\to C_i(\widetilde{M},\C)$ 
as in the previous section, 
we can define the $i$th (combinatorial) Laplace operator 
$\Delta_i:C_i(\widetilde{M},\C)\to C_i(\widetilde{M},\C)$ by 
$$
\Delta_i=\partial_{i+1}\circ\partial_{i+1}^*
+\partial_i^*\circ\partial_i.
$$

Let us suppose that 
all the $L^2$-Betti numbers $b(\Delta_i)$ vanish 
and 
the fundamental group 
$\pi_1M$ is of $\det\geq1$-class. 
Thereby as a generalization of 
the classical Reidemeister torsion, 
the $L^2$-torsion $\tau(M)$ is defined by 

\begin{definition}
{\rm 
$$
\tau(M)=\prod_{i=0}^3\text{det}_{\C\pi_1M}(\Delta_i)^{(-1)^{i+1}i}
\in \mathbb R_{>0}.
$$
}
\end{definition}

As for 
the positivity of Novikov-Shubin invariants 
$\alpha(\Delta_i)$ 
for the Laplace operator $\Delta_i$, 
it is known that 
$\alpha(\Delta_i)>0$ 
holds under some general assumptions 
(see \cite{LL95-1}). 
For example, 
if a compact connected orientable 3-manifold $M$ satisfies 
\begin{enumerate}
 \item $\pi_1M$ is infinite,

 \item $M$ is an irreducible 3-manifold or 
 $S^1\times S^2$ or $\mathbb RP^3\sharp \mathbb RP^3$,
 
 \item if $\partial M\not= \phi$, 
 it consists of tori,
 
 \item if $\partial M= \phi$, 
 $M$ is finitely covered by a 3-manifold which is a hyperbolic, 
Seifert or Haken 3-manifold,
\end{enumerate}
then 
$b(\Delta_i)=0$ and 
$\alpha(\Delta_i)>0$ for each $i$. 
Therefore, 
we see that 
the $L^2$-torsion $\tau(M)$ is 
also well-defined 
in view of these conditions. 

\begin{remark}
The above condition (4) is automatically 
satisfied by Perelman's proof of 
Thurston's Geometrization Conjecture. 
\end{remark}

As a notable property of the $L^2$-torsion, 
it is known that $\log\tau(M)$ can be interpreted 
as Gromov's simplicial volume $||M||$ 
and hyperbolic volume $\mathrm{vol}(M)$ 
(see \cite{Gromov82-1}) of $M$. 
See \cite{LS99-1} for the heart of the proof. 

\begin{theorem}\label{thm:Gromov}
Let $M$ be a compact connected orientable irreducible $3$-manifold 
with an infinite fundamental group such that 
$\partial M$ is empty or a disjoint union of incompressible tori. 
Then 
it holds that 
$$
\log\tau(M)=C||M||,
$$
where $C$ is the universal constant not depending on $M$. 
In particular, 
if $M$ is a hyperbolic $3$-manifold, 
we obtain 
$$
\log\tau(M)=-\frac{1}{3\pi}\mathrm{vol}(M).
$$
\end{theorem}

Next 
we review L\"uck's formula for the $L^2$-torsion 
of $3$-manifolds 
(\cite{Luck94-1} Theorem 2.4). 
From this formula, 
we see that $\log\tau$ is a characteristic polynomial 
of the Magnus representation 
of the mapping class group. 

\begin{theorem}
\label{thm.2.6}
Let $M$ be as in the above theorem. 
We suppose that 
$\partial M$ is non-empty and 
$\pi_1M$ has a deficiency one presentation 
$$
\left\langle s_1,\ldots,s_{n+1}~|~r_1,\ldots,r_n\right\rangle.
$$ 
Put 
$A$ to be the $n\times n$-matrix with entries 
in $\mathbb Z\pi_1M$ obtained from the matrix 
$\left(\partial r_i/\partial s_j\right)$ by deleting 
one of the columns. 
Then the logarithm of the $L^2$-torsion of $M$ is given by 
\begin{align*}
\log\tau (M) 
& = -2\log\mathrm{det}_{\mathbb C\pi_1M}(A) \\
& = -2 n\log K+\sum_{p=1}^\infty\frac1p\text{\rm tr}_{\mathbb C\pi_1M}
\left(\left(I-K^{-2}AA^*\right)^p\right),
\end{align*}
where $K$ is a constant satisfying $K\geq ||R_A||_\infty$.
\end{theorem}

To see a relation between 
the Magnus representation and 
the $L^2$-torsion, 
we describe the above L\"uck's formula 
for a surface bundle over the circle. 

For 
an orientation preserving diffeomorphism 
$\varphi$ of $\Sigma_{g,1}$, 
we form the mapping torus 
$M_\varphi$ by taking the product 
$\Sigma_{g,1}\times [0,1]$ and gluing 
$\Sigma_{g,1}\times\{0\}$ and $\Sigma_{g,1}\times\{1\}$ 
via $\varphi$. 
This gives a surface bundle over $S^1$. 
Its diffeomorphism type is determined by 
the monodromy map $\varphi$, 
and conversely the monodromy map $\varphi$ is determined by a given surface bundle 
up to conjugacy and isotopy. 
Here 
an isotopy fixes setwisely the points 
on the boundary $\partial\Sigma_{g,1}$. 
We take a deficiency one presentation of 
the fundamental group 
$\pi=\pi_1(M_\varphi,*)$, 
$$
\pi
=\left\langle x_1,\ldots,x_{2g},t~|~
r_i:
tx_it^{-1}=\varphi_*(x_i),
~1\leq i\leq 2g\right\rangle,
$$
where 
the base point $*$ of $\pi$ 
and $\Gamma=\pi_1(\Sigma_{g,1},*)$ is the same one 
on the fiber $\Sigma_{g,1}\times\{0\}\subset M_\varphi$ 
and 
$\varphi_*:\Gamma\to\Gamma$ is the automorphism 
induced by $\varphi:\Sigma_{g,1}\to\Sigma_{g,1}$. 
It should be noted that 
$\pi$ is isomorphic to the semi-direct product 
of $\Gamma$ and $\pi_1S^1\cong\Z=\langle t\rangle$. 

Applying the free differential calculus to 
the relations $r_i~(1\leq i\leq 2g)$, 
we obtain the Alexander matrix 
$$
A=\left(\frac{\partial r_i}{\partial x_j}\right)
\in M(2g,\mathbb Z\pi).
$$
Then 
L\"{u}ck's formula for a surface bundle over 
the circle is given by 
\begin{align*}
\log \tau(M_\varphi)
&=-2\log \text{det}_{\mathbb C\pi}(A) \\
&=-4g\log K + 
\sum_{p=1}^\infty\frac{1}{p}\text{tr}_{\C\pi}
\left(\left(I-K^{-2} AA^*\right)^p\right),
\end{align*}
where $K$ is a constant satisfying $K\geq ||R_A||_\infty$.

This formula 
enables us to interpret the $L^2$-torsion 
$\log\tau$ 
of a surface bundle over the circle as 
the characteristic polynomial of 
the Magnus representation 
$r(\varphi)$. 
In fact, 
an easy calculation shows that 
$$
A
=
\left(
\frac{\partial r_i}{\partial x_j}
\right)
=tI-{}^t\overline{r(\varphi)}.
$$
Then 
if we take the Fuglede-Kadison determinant 
in $M(2g,\C\pi)$, 
we have 
\begin{align*}
\text{det}_{\mathbb C\pi}\left(tI-{}^t\overline{r(\varphi)}\right)
&=
\text{det}_{\mathbb C\pi}\left(tI-{}^t\overline{r(\varphi)}\right)^*\\
&=\text{det}_{\mathbb C\pi}\left(t^{-1}I-r(\varphi)\right)
\end{align*}
because $\text{tr}_{\C\pi}(BB^*)=\text{tr}_{\C\pi}(B^*B)$ 
holds. 
Therefore 
the $L^2$-torsion is 
interpreted as the characteristic polynomial 
of $r(\varphi)$. 

\section{Definition of $L^2$-torsion invariants}

As was seen in Section 3, 
L\"uck's formula gives a way to 
calculate the simplicial volume from 
a presentation of the fundamental group. 
However, 
in general, 
it seems to be difficult to 
evaluate the exact values from the formula. 
In this section, 
we introduce a sequence of $L^2$-torsion invariants 
which approximates the original one 
for a surface bundle over the circle. 
See \cite{KMT04-1} for details. 

In order to 
construct such a sequence of $L^2$-torsion invariants, 
we consider the lower central series of $\Gamma$. 
Namely, 
we take the descending infinite sequence 
$$
\Gamma_1=\Gamma\supset \Gamma_2\supset \cdots
\supset\Gamma_k\supset \cdots,
$$
where 
$\Gamma_k=[\Gamma_{k-1},\Gamma_1]$ for $k\geq 2$. 
Let $N_k$ be the $k$th nilpotent quotient $N_k=\Gamma/\Gamma_k$ 
and $p_k:\Gamma\to N_k$ be the natural projection. 

In the previous section, 
we considered a chain complex 
$C_*(\widetilde M_\varphi,\mathbb C)$ 
of $\mathbb C\pi$-modules. 
Instead of this complex, 
we can use the chain complex 
$$
C_*(M_\varphi,l^2(\pi))
=l^2(\pi)\otimes_{\mathbb C\pi}C_*(\widetilde{M}_\varphi,\mathbb C)
$$
to define the same $L^2$-torsion $\tau(M_{\varphi})$. 
This point of view allows us 
to introduce a sequence of 
the $L^2$-torsion invariants. 

The group $\Gamma_k$ is a normal subgroup of $\pi$, 
so that 
we can take the quotient group 
$\pi(k)=\pi/\Gamma_k$. 
It should be noted that 
$\pi(k)$ is isomorphic to the semi-direct product 
$N_k\rtimes\Z$. 
We denote the induced projection 
$\pi\to\pi(k)$ by the same letter $p_k$. 
Thereby 
we can consider the chain complex 
$$
C_*\left(M_\varphi, l^2\left(\pi(k)\right)\right)
= l^2\left(\pi(k)\right)\otimes_{\mathbb C\pi}C_*(\widetilde{M}_\varphi,\mathbb C)
$$
through the projection $p_k$. 
By using the Laplace operator 
$$
\Delta_i^{(k)}:
C_i\left(M_\varphi,l^2\left(\pi(k)\right)\right)
\to
C_i\left(M_\varphi,l^2\left(\pi(k)\right)\right)
$$
on this complex, 
we can formally define 
the $k$th $L^2$-torsion invariant 
$\tau_k(M_\varphi)$ as follows. 

\begin{definition}
{\rm 
$$
\tau_k(M_\varphi)
=\prod_{i=0}^3\text{det}_{\C\pi(k)}(\Delta_i^{(k)})^{(-1)^{i+1}i}.
$$
}
\end{definition}

Of course, 
this definition is well-defined 
if every $L^2$-Betti number $b(\Delta_i^{(k)})$ vanishes 
and 
every $\pi(k)$ is of $\det\geq1$-class. 
The next lemma is easily proved
(see \cite{KMT04-1}, \cite{Luck94-2}). 

\begin{lemma}
The $L^2$-Betti numbers of $\Delta_i^{(k)}$ 
are all zero. 
\end{lemma}

Recall the class $\mathcal{G}$ of groups. 
It is the smallest class of groups 
which contains the trivial group and 
is closed under the following processes: 
(i) amenable quotients, 
(ii) colimits, 
(iii) inverse limits, 
(iv) subgroups and 
(v) quotients with finite kernel 
(see \cite{Luck02-1}, \cite{Schick01-1}). 
It is known that 
$\mathcal{G}$ contains all amenable groups. 
By definition, 
$N_k=\Gamma/\Gamma_k$ is a nilpotent group 
and in particular an amenable group.
Hence 
every $N_k$ belongs to $\mathcal{G}$. 
Further 
for any automorphism 
$\varphi_*:N_k\rightarrow N_k$, 
its mapping torus extension 
({\it HNN}-extension) 
$N_k\rtimes\Z$  also belongs to $\mathcal{G}$. 
Therefore we have 

\begin{lemma}
The group $\pi(k)$ belongs to $\mathcal{G}$. 
\end{lemma}

As a result, 
we can conclude that 
our $L^2$-torsion invariants $\tau_k$ 
can be defined for any $k\geq 1$ and 
they are all homotopy invariants 
(see \cite{Luck02-1}, \cite{Schick01-1}). 

Now 
let us describe a formula 
of the $k$th $L^2$-torsion invariant 
$\tau_k(M_\varphi)$ 
and 
establish a relation to 
the Magnus representation of 
the mapping class group. 
Let 
${p_k}_*:\C\pi\to\C\pi(k)$ be 
an induced homomorphism over the group rings. 
For 
$k\geq 1$, 
we put 
$$
A_k=
\left({p_k}_*\left(\frac{\partial r_i}{\partial x_j}\right)\right)
\in M(2g,\mathbb C\pi(k)).
$$
Moreover 
we fix a constant $K_k$ satisfying $K_k\geq ||R_{A_k}||_\infty$. 
Then 
we have 
\begin{align*}
\log \tau_k(M_\varphi)
&=-2\log \text{det}_{\mathbb C\pi(k)}(R_{A_k}) \\
&=-4g\log K_k + 
\sum_{p=1}^\infty\frac{1}{p}\text{tr}_{\C\pi(k)}
\left(\left(I-K_k^{-2} A_k{A_k}^*\right)^p\right),
\end{align*}
by virtue of the same argument as 
Theorem \ref{thm.2.6}.

For the $k$th invariant $\tau_k$, 
we have taken the lower central series 
$\{\Gamma_k\}$ of $\Gamma$ 
and the nilpotent quotients 
$\{N_k
\}$. 
These quotients induce a sequence of representations 
(more precisely, 
crossed homomorphisms) 
$$
r_k:\M_{g,1}\to GL(2g,\Z N_k)
$$
for $k\geq 1$ (see \cite{Morita93-1}). 
They 
naively approximate the original 
Magnus representation 
$r:\M_{g,1}\to GL(2g,\Z\Gamma)$. 
By the similar observation as before, 
the $k$th invariant $\log\tau_k(M_\varphi)$ 
can be regarded as the characteristic polynomial of $r_k(\varphi)$ 
with respect to the Fuglede-Kadison determinant 
in $M(2g,\C\pi(k))$. 

From 
the viewpoint of the Magnus representation 
of the mapping class group, 
it seems natural to raise 
the following problem. 

\begin{problem}\label{conj.3.7}
Show that 
the sequence $\left\{ \tau_k(M_\varphi)\right\}$ 
converges to $\tau(M_\varphi)$ 
when we take the limit on $k$. 
\end{problem}

In general, 
such an approximation problem for 
the $L^2$-torsion seems to be difficult. 
However, 
similar convergence results are known 
for the $L^2$-Betti numbers. 
In fact, 
L\"uck shows in \cite{Luck94-3} 
a theorem relating $L^2$-Betti numbers to 
ordinary Betti numbers of finite coverings. 
This result is generalized to 
more general settings by 
Schick in \cite{Schick01-1}. 

As for 
the Fuglede-Kadison determinant, 
L\"uck proves in \cite{Luck02-1} 
the following. 
Let 
$f:\Q[\Z]\to\Q[\Z]$ 
be the $\Q[\Z]$-map 
given by multiplication 
with $p(t)\in \Q[\Z]$ and 
$f_{(2)}:l^2(\Z)\to l^2(\Z)$ 
be the linear operator obtained from $f$ 
by tensoring with $l^2(\Z)$ 
over $\Q[\Z]$. 
Further 
let $f_{[n]}:\C[\Z/n]\to\C[\Z/n]$ 
be the linear operator obtained from $f$ 
by taking the tensor product with 
$\C[\Z/n]$ over $\Q[\Z]$. 
We then 
get an approximation result: 
$$
\log\mathrm{det}_{\C[\Z]}\left(f_{(2)}\right)
=
\lim_{n\to\infty}
\frac{\log\mathrm{det}_{\C[\Z/n]}\left(f_{[n]}\right)}{n}
$$
(see \cite{KMT03-2} for a similar statement). 
In \cite{Luck02-1} 
L\"uck also points out 
that 
there exists a purely algebraic example 
where 
Fuglede-Kadison determinants 
do not converge. 

On the other hand, 
in general, 
we have at least an inequality 
for the Fuglede-Kadison determinant 
in the limit statement 
(see \cite{Schick01-1}). 
That is, 
for the operator $R_{A_k}$ 
we see that 
$$
\log\mathrm{det}_{\C\pi}(R_{A})
\geq
\limsup_k \log\mathrm{det}_{\C\pi(k)}(R_{A_k})
$$
holds. 
In the last section, 
we shall discuss Problem \ref{conj.3.7} again and give 
an affirmative answer under certain conditions. 

\section{Formulas of $\tau_1$ and $\tau_2$}

In this section, 
we give explicit formulas of 
the first two invariants 
of a sequence of our $L^2$-torsion invariants. 
They are really computable formulas, 
so that we can make a systematic calculation for low genus cases. 
In particular, 
we compare them 
with hyperbolic volumes. 
The results discussed here 
are a summary of our previous paper \cite{KMT04-1} 
(see also \cite{KMT03-1}, \cite{KMT03-2}). 

First 
we consider the Magnus representation 
$$
r_1:\M_{g,1}\to GL(2g,\Z N_1).
$$
Here 
$N_1=\Gamma/\Gamma_1$ is the trivial group 
and then 
the above representation is 
the same as 
the usual homological action 
of $\M_{g,1}$ on $H_1(\Sigma_{g,1},\Z)$. 
Namely 
we have the representation 
$$
r_1:\M_{g,1}\to
\mathrm{Aut}\left(H_1(\Sigma_{g,1},\Z),\langle\ ,\ \rangle\right)
\cong \mathrm{Sp}(2g,\Z),
$$
where 
$\langle\ ,\ \rangle$ denotes the intersection form 
on the first homology group. 
Further 
$\pi(1)=\pi/\Gamma_1\cong
\Z=\langle t\rangle$ 
and its group ring $\C\langle t \rangle$ is  
a commutative Laurent polynomial ring 
$\C[t, t^{-1}]$. 
Then the matrix 
$A_1$ is nothing but the usual characteristic matrix 
of ${}^tr_1(\varphi)$. 
In this case, 
it is described by the usual determinant 
for a matrix with commutative entries. 

In order to state the theorem, 
we recall a definition from number theory 
(see \cite{Everest98-1} and its references). 
For 
a Laurent polynomial 
$F(\mathbf{t})\in \C[t_1^{\pm1},\ldots,t_n^{\pm1}]$, 
the Mahler measure of $F$ is defined by 
$$
m(F)
=\int_{0}^{1}\cdots\int_0^1\log 
\left|F(e^{2\pi\sqrt{-1}\theta_1},\ldots,
e^{2\pi\sqrt{-1}\theta_n})\right|
d{\theta_1}\cdots d{\theta_n}, 
$$
where 
we assume that undefined terms are omitted. 
Namely 
we define the integrand to be zero 
whenever 
we hit a zero of $F$. 

\begin{theorem}[\cite{KMT04-1}]\label{thm.5.1}
The logarithm of the first invariant $\tau_1$ is given by 
$$
\log\tau_1(M_{\varphi})
= -2 m \left(\Delta_{r_1(\varphi)}\right),
$$
where 
$\Delta_{r_1(\varphi)}(t)
=\det A_1=\mathrm{det}(t I-r_1(\varphi))$. 
Moreover 
if 
$\Delta_{r_1(\varphi)}(t)$ 
has a factorization 
$\Delta_{r_1(\varphi)}(t)
= \prod_{i=1}^{2g}(t-\alpha_i)\ (\alpha_i\in\C)$, 
then we have 
$$
\log\tau_1(M_\varphi)
=-2\sum_{i=1}^{2g} \log \mathrm{max} \{ 1, |\alpha_i|\} .
$$
\end{theorem}

\begin{remark}
{\rm 
In other words, 
$\log\tau_1(M_\varphi)$ 
is given by the integral of 
the Alexander polynomial of 
$M_\varphi$ over the circle $S^1$ 
(see \cite{Luck94-1}, 
for the exterior of a knot $K$ in the $3$-sphere $S^3$). 
Further, 
$\log\tau_1(M_\varphi)$ 
can be described by the asymptotic behavior 
of the order of the first homology group 
of a cyclic covering 
(see \cite{KMT03-2}). 
}
\end{remark}

The point of the proof is 
to identify the Hilbert space $l^2(\Z)$ 
with $L^2(\R/\Z)$ 
in terms of the Fourier transforms. 
Then 
the $\C\langle t\rangle$-trace 
$\text{tr}_{\C\langle t\rangle}:l^2(\Z)\to\C$ 
can be realized as the integration 
$$
L^2(\R/\Z)\ni f(\theta)\mapsto \int_0^1 f(\theta)d\theta \in \C
$$
(see \cite{KMT04-1} for details). 
From this description 
and Kronecker's theorem 
(\cite{Everest98-1} Theorem 2), 
we obtain 
a certain vanishing theorem of 
the first invariant. 

\begin{corollary}\label{cor.5.3}
The logarithm of 
$\tau_1(M_\varphi)$ vanishes if and only if 
every eigenvalue of $r_1(\varphi)\in \mathrm{Sp}(2g,\mathbb Z)$
is a root of unity.
\end{corollary}

This corollary seems to be interesting 
from the viewpoint of Problem \ref{conj.3.7}. 
Because in some case, 
we can say that 
the first invariant $\tau_1$ already approximates 
the simplicial volume. 
In particular, 
Corollary \ref{cor.5.3} implies that 
a torus bundle $M_\varphi~(g=1)$ with 
a hyperbolic structure (namely, 
$|\text{tr}\left(r_1(\varphi)\right)|\geq 3$) 
has always non-trivial $L^2$-torsion invariant 
$\tau_1(M_\varphi)$. 
Summing up, 
we have 

\begin{corollary}
For any $\varphi\in \M_{1,1}$, 
its mapping torus $M_\varphi$ 
admits a hyperbolic structure 
if and only if $M_\varphi$ has a non-trivial 
$L^2$-torsion invariant $\tau_1(M_\varphi)$.
\end{corollary}

Therefore, 
the first invariant $\tau_1$ already approximates 
the simplicial volume in genus one case. 

\begin{remark}
It is known that 
if the characteristic polynomial of 
$r_1(\varphi)\in \mathrm{Sp}(2g,\Z)$ is 
irreducible over $\Z$, 
has no roots of unity as eigenvalues and 
is not a polynomial in $t^n$ for any $n>1$, 
then 
$\varphi$ is pseudo-Anosov 
(see Casson-Bleiler \cite{CB88-1}). 
In this case, 
$\mathrm{vol}(M_\varphi)\not=0$ and 
further $\log\tau_1(M_\varphi)\not=0$ by 
Corollary \ref{cor.5.3}.
\end{remark}

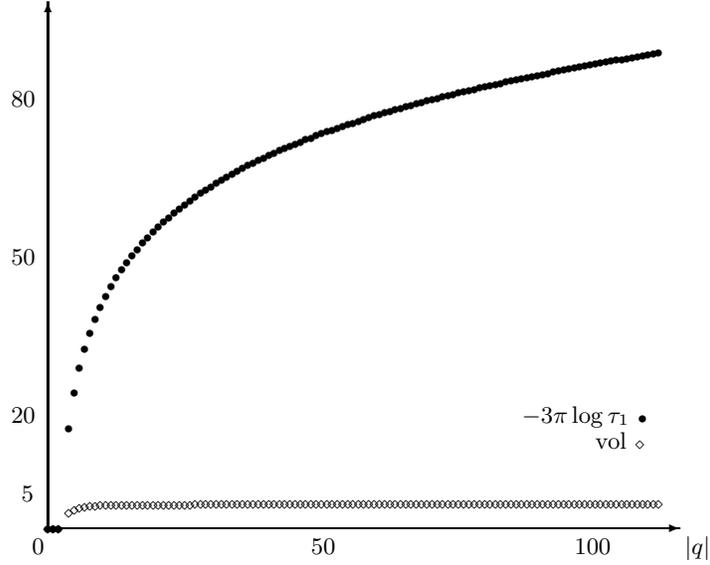
\begin{figure}
\setlength{\unitlength}{0.7mm}
\begin{picture}(115,100)(1,0)

\put(90,20){{\small $-3\pi\log\tau_1$ {\tiny$\bullet$}}}
\put(104,15){{\small vol {\tiny $\diamond$}}}

\put(0,0){\vector(1,0){120}}
\put(0,0){\vector(0,1){100}}

\put(-3,-5){{\small 0}}
\put(50,-5){{\small 50}}
\put(100,-5){{\small 100}}
\put(121,-5){{\small $|q|$}}

\put(-7,80){{\small 80}}
\put(-7,50){{\small 50}}
\put(-7,20){{\small 20}}
\put(-5,5){{\small 5}}

\put(-1,-1){{\tiny$\bullet$}}
\put(0,-1){{\tiny$\bullet$}}
\put(1,-1){{\tiny$\bullet$}}
\put(3,18.14){{\tiny$\bullet$}}
\put(4,24.82){{\tiny$\bullet$}}
\put(5,29.53){{\tiny$\bullet$}}
\put(6,33.22){{\tiny$\bullet$}}
\put(7,36.28){{\tiny$\bullet$}}
\put(8,38.89){{\tiny$\bullet$}}
\put(9,41.17){{\tiny$\bullet$}}
\put(10,43.21){{\tiny$\bullet$}}
\put(11,45.04){{\tiny$\bullet$}}
\put(12,46.70){{\tiny$\bullet$}}
\put(13,48.23){{\tiny$\bullet$}}
\put(14,49.64){{\tiny$\bullet$}}
\put(15,50.96){{\tiny$\bullet$}}
\put(16,52.1879998735){{\tiny$\bullet$}}
\put(17,53.339248833){{\tiny$\bullet$}}
\put(18,54.4237752394){{\tiny$\bullet$}}
\put(19,55.4489340545){{\tiny$\bullet$}}
\put(20,56.4209214237){{\tiny$\bullet$}}
\put(21,57.3450070665){{\tiny$\bullet$}}
\put(22,58.2257109613){{\tiny$\bullet$}}
\put(23,59.0669396703){{\tiny$\bullet$}}
\put(24,59.8720929219){{\tiny$\bullet$}}
\put(25,60.6441479351){{\tiny$\bullet$}}
\put(26,61.3857268559){{\tiny$\bullet$}}
\put(27,62.0991512184){{\tiny$\bullet$}}
\put(28,62.7864863245){{\tiny$\bullet$}}
\put(29,63.4495777046){{\tiny$\bullet$}}
\put(30,64.0900813018){{\tiny$\bullet$}}
\put(31,64.7094886357){{\tiny$\bullet$}}
\put(32,65.3091479165){{\tiny$\bullet$}}
\put(33,65.8902818691){{\tiny$\bullet$}}
\put(34,66.4540028652){{\tiny$\bullet$}}
\put(35,67.001325837){{\tiny$\bullet$}}
\put(36,67.5331793532){{\tiny$\bullet$}}
\put(37,68.0504151613){{\tiny$\bullet$}}
\put(38,68.5538164472){{\tiny$\bullet$}}
\put(39,69.0441050106){{\tiny$\bullet$}}
\put(40,69.5219475257){{\tiny$\bullet$}}
\put(41,69.9879610217){{\tiny$\bullet$}}
\put(42,70.4427176988){{\tiny$\bullet$}}
\put(43,70.886749173){{\tiny$\bullet$}}
\put(44,71.3205502291){{\tiny$\bullet$}}
\put(45,71.74458215){{\tiny$\bullet$}}
\put(46,72.1592756775){{\tiny$\bullet$}}
\put(47,72.5650336525){{\tiny$\bullet$}}
\put(48,72.9622333763){{\tiny$\bullet$}}
\put(49,73.3512287261){{\tiny$\bullet$}}
\put(50,73.7323520569){{\tiny$\bullet$}}
\put(51,74.1059159129){{\tiny$\bullet$}}
\put(52,74.4722145732){{\tiny$\bullet$}}
\put(53,74.8315254481){{\tiny$\bullet$}}
\put(54,75.1841103457){{\tiny$\bullet$}}
\put(55,75.530216621){{\tiny$\bullet$}}
\put(56,75.8700782211){{\tiny$\bullet$}}
\put(57,76.2039166379){{\tiny$\bullet$}}
\put(58,76.5319417771){{\tiny$\bullet$}}
\put(59,76.8543527532){{\tiny$\bullet$}}
\put(60,77.1713386169){{\tiny$\bullet$}}
\put(61,77.4830790227){{\tiny$\bullet$}}
\put(62,77.7897448414){{\tiny$\bullet$}}
\put(63,78.0914987245){{\tiny$\bullet$}}
\put(64,78.3884956224){{\tiny$\bullet$}}
\put(65,78.6808832639){{\tiny$\bullet$}}
\put(66,78.9688025983){{\tiny$\bullet$}}
\put(67,79.2523882037){{\tiny$\bullet$}}
\put(68,79.531768666){{\tiny$\bullet$}}
\put(69,79.8070669292){{\tiny$\bullet$}}
\put(70,80.0784006211){{\tiny$\bullet$}}
\put(71,80.3458823554){{\tiny$\bullet$}}
\put(72,80.609620013){{\tiny$\bullet$}}
\put(73,80.8697170033){{\tiny$\bullet$}}
\put(74,81.1262725085){{\tiny$\bullet$}}
\put(75,81.3793817106){{\tiny$\bullet$}}
\put(76,81.6291360039){{\tiny$\bullet$}}
\put(77,81.8756231935){{\tiny$\bullet$}}
\put(78,82.1189276808){{\tiny$\bullet$}}
\put(79,82.359130637){{\tiny$\bullet$}}
\put(80,82.5963101663){{\tiny$\bullet$}}
\put(81,82.8305414583){{\tiny$\bullet$}}
\put(82,83.0618969315){{\tiny$\bullet$}}
\put(83,83.2904463673){{\tiny$\bullet$}}
\put(84,83.5162570375){{\tiny$\bullet$}}
\put(85,83.7393938222){{\tiny$\bullet$}}
\put(86,83.9599193227){{\tiny$\bullet$}}
\put(87,84.1778939666){{\tiny$\bullet$}}
\put(88,84.3933761072){{\tiny$\bullet$}}
\put(89,84.606422118){{\tiny$\bullet$}}
\put(90,84.8170864804){{\tiny$\bullet$}}
\put(91,85.0254218684){{\tiny$\bullet$}}
\put(92,85.2314792268){{\tiny$\bullet$}}
\put(93,85.4353078467){{\tiny$\bullet$}}
\put(94,85.6369554362){{\tiny$\bullet$}}
\put(95,85.8364681874){{\tiny$\bullet$}}
\put(96,86.0338908398){{\tiny$\bullet$}}
\put(97,86.2292667409){{\tiny$\bullet$}}
\put(98,86.4226379034){{\tiny$\bullet$}}
\put(99, 86.6140450591){{\tiny$\bullet$}}
\put(100,86.8035277106){{\tiny$\bullet$}}
\put(101,86.9911241805){{\tiny$\bullet$}}
\put(102,87.1768716579){{\tiny$\bullet$}}
\put(103,87.3608062421){{\tiny$\bullet$}}
\put(104,87.5429629857){{\tiny$\bullet$}}
\put(105,87.723375934){{\tiny$\bullet$}}
\put(106,87.9020781634){{\tiny$\bullet$}}
\put(107,88.079101818){{\tiny$\bullet$}}
\put(108,88.2544781442){{\tiny$\bullet$}}
\put(109,88.4282375237){{\tiny$\bullet$}}
\put(110,88.6004095054){{\tiny$\bullet$}}
\put(111,88.7710228353){{\tiny$\bullet$}}
\put(112,88.9401054854){{\tiny$\bullet$}}
\put(113,89.1076846815){{\tiny$\bullet$}}
\put(114,89.273786929){{\tiny$\bullet$}}
\put(115,89.4384380386){{\tiny$\bullet$}}

\put(-1,-1){{\tiny $\diamond$}}
\put(0,-1){{\tiny $\diamond$}}
\put(1,-1){{\tiny $\diamond$}}
\put(3	,2.02988321281930){{\tiny $\diamond$}} 
\put(4,2.66674478344910 ){{\tiny $\diamond$}}
\put(5	,2.98912028295000 ){{\tiny $\diamond$}}
\put(6	,3.17729327860033 ){{\tiny $\diamond$}}
\put(7	,3.29690241432664 ){{\tiny $\diamond$}}
\put(8	,3.37759740823144 ){{\tiny $\diamond$}}
\put(9	,3.43454088590256 ){{\tiny $\diamond$}}
\put(10,3.47617398923899){{\tiny $\diamond$}} 
\put(11,3.50750355950584 ){{\tiny $\diamond$}}
\put(12,3.53165135729638 ){{\tiny $\diamond$}}
\put(13,3.55064458472050 ){{\tiny $\diamond$}}
\put(14,3.56584510355300 ){{\tiny $\diamond$}}
\put(15,3.57819478022670 ){{\tiny $\diamond$}}
\put(16,3.58836116484430 ){{\tiny $\diamond$}}
\put(17,3.59682804103340 ){{\tiny $\diamond$}}
\put(18,3.60395268263800 ){{\tiny $\diamond$}}
\put(19,3.61000348403880 ){{\tiny $\diamond$}}
\put(20,3.61518516184410 ){{\tiny $\diamond$}}
\put(21,3.61965601470680 ){{\tiny $\diamond$}}
\put(22,3.62353997319740 ){{\tiny $\diamond$}}
\put(23,3.62693515911146 ){{\tiny $\diamond$}}
\put(24,3.62992006280211 ){{\tiny $\diamond$}}
\put(25,3.63255806429380 ){{\tiny $\diamond$}}
\put(26,3.63490078747160 ){{\tiny $\diamond$}}
\put(27,3.63699061798150 ){{\tiny $\diamond$}}
\put(28,3.63886261479507 ){{\tiny $\diamond$}}
\put(29,3.64054597647090 ){{\tiny $\diamond$}}
\put(30,3.64206517667010 ){{\tiny $\diamond$}}
\put(31,3.64344085140820 ){{\tiny $\diamond$}}
\put(32,3.64469049811070 ){{\tiny $\diamond$}}
\put(33,3.64582903069350 ){{\tiny $\diamond$}}
\put(34,3.64686922354376 ){{\tiny $\diamond$}}
\put(35,3.64782206907770 ){{\tiny $\diamond$}}
\put(36,3.64869706755880 ){{\tiny $\diamond$}}
\put(37,3.64950246343950 ){{\tiny $\diamond$}}
\put(38,3.65024543920270 ){{\tiny $\diamond$}}
\put(39,3.65093227520780 ){{\tiny $\diamond$}}
\put(40,3.65156848218390 ){{\tiny $\diamond$}}
\put(41,3.65215891158500 ){{\tiny $\diamond$}}
\put(42,3.65270784793340 ){{\tiny $\diamond$}}
\put(43,3.65321908642980 ){{\tiny $\diamond$}}
\put(44,3.65369599845720 ){{\tiny $\diamond$}}
\put(45,3.65414158708500 ){{\tiny $\diamond$}}
\put(46,3.65455853428150 ){{\tiny $\diamond$}}
\put(47,3.65494924121880 ){{\tiny $\diamond$}}
\put(48,3.65531586280044 ){{\tiny $\diamond$}}
\put(49,3.65566033733950 ){{\tiny $\diamond$}}
\put(50,3.65598441214847 ){{\tiny $\diamond$}}
\put(51,3.65628966567280 ){{\tiny $\diamond$}}
\put(52,3.65657752669060 ){{\tiny $\diamond$}}
\put(53,3.65684929101423 ){{\tiny $\diamond$}}
\put(54,3.65710613605900 ){{\tiny $\diamond$}}
\put(55,3.65734913358391 ){{\tiny $\diamond$}}
\put(56,3.65757926086160 ){{\tiny $\diamond$}}
\put(57,3.65779741049570 ){{\tiny $\diamond$}}
\put(58,3.65800439906890 ){{\tiny $\diamond$}}
\put(59,3.65820097477790 ){{\tiny $\diamond$}}
\put(60,3.65838782418960 ){{\tiny $\diamond$}}
\put(61,3.65856557823077 ){{\tiny $\diamond$}}
\put(62,3.65873481751060 ){{\tiny $\diamond$}}
\put(63,3.65889607705800 ){{\tiny $\diamond$}}
\put(64,3.65904985054720 ){{\tiny $\diamond$}}
\put(65,3.65919659407390 ){{\tiny $\diamond$}}
\put(66,3.65933672953450 ){{\tiny $\diamond$}}
\put(67,3.65947064765657 ){{\tiny $\diamond$}}
\put(68,3.65959871072100 ){{\tiny $\diamond$}}
\put(69,3.65972125501066 ){{\tiny $\diamond$}}
\put(70,3.65983859301620 ){{\tiny $\diamond$}}
\put(71,3.65995101542760 ){{\tiny $\diamond$}}
\put(72,3.66005879293330 ){{\tiny $\diamond$}}
\put(73,3.66016217784870 ){{\tiny $\diamond$}}
\put(74,3.66026140559200 ){{\tiny $\diamond$}}
\put(75,3.66035669602425 ){{\tiny $\diamond$}}
\put(76,3.66044825466550 ){{\tiny $\diamond$}}
\put(77,3.66053627380300 ){{\tiny $\diamond$}}
\put(78,3.66062093350200 ){{\tiny $\diamond$}}
\put(79,3.66070240252210 ){{\tiny $\diamond$}}
\put(80,3.66078083916200 ){{\tiny $\diamond$}}
\put(81,3.66085639202580 ){{\tiny $\diamond$}}
\put(82,3.66092920072800 ){{\tiny $\diamond$}}
\put(83,3.66099939653620 ){{\tiny $\diamond$}}
\put(84,3.66106710296500 ){{\tiny $\diamond$}}
\put(85,3.66113243631723 ){{\tiny $\diamond$}}
\put(86,3.66119550618390 ){{\tiny $\diamond$}}
\put(87,3.66125641590400 ){{\tiny $\diamond$}}
\put(88,3.66131526298530 ){{\tiny $\diamond$}}
\put(89,3.66137213949920 ){{\tiny $\diamond$}}
\put(90,3.66142713243580 ){{\tiny $\diamond$}}
\put(91,3.66148032403910 ){{\tiny $\diamond$}}
\put(92,3.66153179211370 ){{\tiny $\diamond$}}
\put(93,3.66158161030917 ){{\tiny $\diamond$}}
\put(94,3.66162984838350 ){{\tiny $\diamond$}}
\put(95,3.66167657244746 ){{\tiny $\diamond$}}
\put(96,3.66172184519056 ){{\tiny $\diamond$}}
\put(97,3.66176572609170 ){{\tiny $\diamond$}}
\put(98,3.66180827161410 ){{\tiny $\diamond$}}
\put(99,3.66184953538740 ){{\tiny $\diamond$}}
\put(100,3.66188956837600 ){{\tiny $\diamond$}}
\put(101,3.66192841903680 ){{\tiny $\diamond$}}
\put(102,3.66196613346550 ){{\tiny $\diamond$}}
\put(103,3.66200275553360 ){{\tiny $\diamond$}}
\put(104,3.66203832701580 ){{\tiny $\diamond$}}
\put(105,3.66207288770900 ){{\tiny $\diamond$}}
\put(106,3.66210647554380 ){{\tiny $\diamond$}}
\put(107,3.66213912668820 ){{\tiny $\diamond$}}
\put(108,3.66217087564543 ){{\tiny $\diamond$}}
\put(109,3.66220175534470 ){{\tiny $\diamond$}}
\put(110,3.66223179226600 ){{\tiny $\diamond$}}
\put(111,3.66226103132330 ){{\tiny $\diamond$}}
\put(112,3.66228948633300 ){{\tiny $\diamond$}}
\put(113,3.66231718969260 ){{\tiny $\diamond$}}
\put(114,3.66234416763930 ){{\tiny $\diamond$}}
\put(115,3.66237044527730 ){{\tiny $\diamond$}}
\end{picture}
\caption{$-3\pi\log\tau_1(M_\varphi)$ and $\mathrm{vol}(M_\varphi)$ vs. $|q|$}              
\end{figure}

\begin{example}
{\rm 
It is well-known 
that the mapping class group of 
the two dimensional torus $T^2=\R^2/\Z^2$ 
is isomorphic to $SL(2,\Z)$. 
Taking a matrix 
$\begin{pmatrix}
q&1\\
-1&0
\end{pmatrix}\in SL(2,\Z)$, 
it gives a diffeomorphism $\varphi$ on $T^2$. 
We may assume that 
it is the identity on some embedded disk 
by an isotopic deformation and it gives an element of $\M_{1,1}$. 
We use the same symbol $\varphi$ for this mapping class. 
An easy calculation shows that 
$$
r_1(\varphi)
=\begin{pmatrix}
q&1\\
-1&0
\end{pmatrix}$$
and 
$$
\Delta_{r_1(\varphi)}(t)
=\det(tI-r_1(\varphi))
=t^2-qt+1.
$$
We put 
$\xi_\pm=(q\pm\sqrt{q^2-4})/2$ 
(the eigenvalues of the matrix $r_1(\varphi)$). 
If 
$|q|\leq2$, 
then $|\xi_\pm|=1$. 
Hence 
$\log\tau_1(M_\varphi)=0$ 
in these cases. 
On the other hand, 
either $|\xi_+|$ or $|\xi_-|$ 
is greater than one 
when $|q|\geq3$, 
so that 
$M_\varphi$ has a non-trivial 
$L^2$-torsion invariant $\tau_1$ 
in these cases. 
In fact, 
the logarithm of the first invariant 
is given by 
$$
\log\tau_1(M_\varphi)
=
-2\log\max\left\{|\xi_+|,|\xi_-|\right\}.
$$

The values of $\log\tau_1$ 
for the traces $q$ and $-q$ are 
the same, 
so that 
it is a function of $|\mathrm{tr}(r_1(\varphi))|$. 
We put a graph of the $L^2$-torsion invariant 
$-3\pi\log\tau_1(M_\varphi)$ and the hyperbolic 
volume $\mathrm{vol}(M_\varphi)$ 
as a function of $|q|$ 
in Fig 1. 
}
\end{example}

\begin{example}
{\rm 
Next 
we consider the genus two case. 
Let 
$t_1,\ldots, t_5$ be 
the Lickorish-Humphries 
generators of $\M_{2,1}$. 
We take the element 
$\varphi=t_1 t_3 {t_5}^2 t_2^{-1} t_4^{-1}
\in \M_{2,1}$. 
As was shown in \cite{CB88-1}, 
the characteristic polynomial of $r(\varphi)$ 
is 
\begin{align*}
\Delta_{r_1(\varphi)}(t)
&=\det(tI-r_1(\varphi))\\
&=t^4-9t^3+21t^2-9t+1
\end{align*}
and 
irreducible over $\Z$. 
Moreover 
it has no roots of unity as zeros. 
Hence, 
$\varphi$ is pseudo-Anosov and 
$M_\varphi$ has a non-trivial $L^2$-torsion 
invariant $\tau_1(M_\varphi)$. 
In fact, 
we have 
$$
-3\pi\log\tau_1(M_\varphi)
=52.954....\quad
\mathrm{and}\quad
\mathrm{vol}(M_\varphi)=11.466....
$$
}
\end{example}

\begin{remark}
In the above two examples, 
we used SnapPea \cite{Weeks-1}
to compute the hyperbolic volumes. 
\end{remark}

Now 
in the following, 
we consider the second invariant $\tau_2$. 
In the case of genus one, 
we can prove the vanishing of 
$\log\tau_2(M_\varphi)$. 

\begin{theorem}[\cite{KMT03-2}]
$\log\tau_2(M_\varphi)=0$ for any 
$\varphi\in\M_{1,1}$. 
\end{theorem}

This follows from the fact that 
the group $\pi(2)$ is isomorphic to 
the fundamental group of a closed 
torus bundle over the circle. 
Such a $3$-manifold admits 
no hyperbolic structures, 
so that 
the original $L^2$-torsion is trivial 
and 
we obtain the assertion.

On the other hand, 
in the case of 
$g\geq2$, 
it is difficult to describe $\log\tau_2$ explicitly 
on the full mapping class group $\M_{g,1}$. 
However, 
we can do it on the Torelli group. 
Let 
$\varphi$ be an element of the Torelli group $\I_{g,1}$, 
namely 
$\varphi$ acts trivially on the first homology group 
$H_1(\Sigma_{g,1},\Z)$. 
Then 
we notice that 
$\log\tau_1(M_\varphi)=0$ holds 
for any $\varphi\in \I_{g,1}$ 
(see Corollary \ref{cor.5.3}). 
To give an explicit formula of 
$\log\tau_2$, 
we consider the Magnus representation 
$$
r_2:\M_{g,1}\to GL(2g,\Z N_2),
$$
where 
$N_2=\Gamma/[\Gamma,\Gamma]\cong H_1(\Sigma_{g,1},\Z).$ 
If we restrict $r_2$ to the Torelli group $\I_{g,1}$, 
this is really a homomorphism (see \cite{Morita93-1} Corollary 5.4). 
Then 
our formula for the second $L^2$-torsion invariant is the following. 
The proof is similar to 
one for Theorem \ref{thm.5.1}.

\begin{theorem}[\cite{KMT04-1}]\label{thm:Torelli}
For any mapping class $\varphi\in \I_{g,1}$, 
the logarithm of the second $L^2$-torsion invariant 
$\tau_2(M_\varphi)$ 
is given by
$$
\log\tau_2(M_\varphi)
=-2 m\left(\Delta_{r_2(\varphi)}\right),
$$
where
$\Delta_{r_2(\varphi)}(y_1,\ldots,y_{2g},t)
=\det A_2
=\det(tI-\overline{r_2(\varphi)})$ 
and $y_i$ denotes the homology class 
corresponding to $x_i$. 
\end{theorem}

Now 
we suppose 
$F(\mathbf{t})\in \Z[t_1^{\pm1},\ldots,t_n^{\pm1}]$ 
is primitive. 
We 
define $F$ to be a generalized cyclotomic polynomial 
if it is a monomial times a product of one-variable 
cyclotomic polynomials evaluated at monomials. 

The next corollary immediately follows from 
the theorem of Boyd, Lawton and Smyth 
(see \cite{Everest98-1} Theorem 4). 

\begin{corollary}\label{cor.6.2}
For any mapping class $\varphi\in \I_{g,1}$, 
$\log\tau_2(M_\varphi)=0$ 
if and only if 
$\Delta_{r_2(\varphi)}$ is a generalized cyclotomic polynomial.
\end{corollary}

As a typical element of the Torelli group $\I_{g,1}$, 
we first consider a BSCC-map 
$\varphi_h\ (1\leq h \leq g)$ of genus $h$. 
That is, 
a Dehn twist along a bounding simple closed curve on 
$\Sigma_{g,1}$ which separates $\Sigma_{g,1}$ into 
$\Sigma_{h,1}$ and genus $g-h$ surface with two boundaries. 
We then see from \cite{Su0} Corollary 4.3 that 
$\Delta_{r_2(\varphi_h)}
=(t-1)^{2g}$. 
This is clearly 
a generalized cyclotomic polynomial, 
so that $\log\tau_2(M_{\varphi_h})=0$. 

Second 
we consider a BP-map $\psi_h=D_c D_{c'}^{-1}$ of 
genus $h\ (1\leq h \leq g-1)$, 
where $c$ and $c'$ are disjoint homologous simple closed curves 
on $\Sigma_{g,1}$ 
and 
$D_c$ denotes the Dehn twist along $c$. 
Since 
$$
\Delta_{r_2(\psi_h)}
=(t-1)^{2g-2h}(t-y_{g+h+1})^{2h}
$$ 
holds (see \cite{Su0}), 
where $y_{g+h+1}$ denotes the homology class 
corresponding to the $(h+1)$th meridian of $\Sigma_{g,1}$, 
we also have $\log\tau_2(M_{\psi_h})=0$. 

The next example shows the non-triviality of the second 
$L^2$-torsion invariant $\log\tau_2$. 

\begin{example}
{\rm 
Let 
$\varphi=t_3\varphi_1t_3^{-1}\varphi_1\in \I_{2,1}$. 
Then 
we see from a computation in \cite{Su0} that 
$$
\Delta_{r_2(\varphi)}
=(t-1)^4+t(t-1)^2(y_{1}-2+{y_{1}}^{-1})(y_{2}-2+{y_{2}}^{-1}).
$$
This is not a generalized cyclotomic polynomial, 
so that the mapping torus $M_\varphi$ has 
a non-trivial $L^2$-torsion invariant $\tau_2(M_\varphi)$. 
In fact 
we can numerically compute it 
by means of Lawton's result (see \cite{Lawton83-1}). 
More precisely 
we have 
\begin{align*}
-3\pi\log\tau_2(M_\varphi)
&=6\pi m\left(\Delta_{r_2(\varphi)}\right)\\
&=6\pi\lim_{r\to\infty}m\left(\Delta_{r_2(\varphi)}(u,u,u^r)\right)\\
&=19.28....
\end{align*}
}
\end{example}

\section{Vanishing of $\log\tau_k$ for reducible mapping classes}

From the Nielsen-Thurston theory (see \cite{CB88-1}), 
the mapping classes of a surface are 
classified into the following three types: 
(i) periodic, 
(ii) reducible and 
(iii) pseudo-Anosov. 
In our point of view, 
the most interesting object is a pseudo-Anosov map $\varphi$. 
Because 
the corresponding mapping torus $M_\varphi$ has 
non-trivial hyperbolic volume. 

In this final section, 
we show 
two vanishing theorems for $\log\tau_k$.  
We introduced an infinite sequence $\{\tau_k\}$ 
as an approximation of the hyperbolic volume. 
Thus 
if it behaves well with the index $k$, 
we ought to prove 
$$
\lim_{k\to\infty}
\log\tau_k=0
$$ 
for non-hyperbolic 3-manifolds 
(see Problem \ref{conj.3.7}). 
As a first step of this observation, 
we obtain the following. 

\begin{theorem}
If $\varphi\in \M_{g,1}$ is the product 
of Dehn twists along any disjoint non-separating 
simple closed curves on $\Sigma_{g,1}$ 
which are mutually non-homologous, 
then 
$\log \tau_k(M_\varphi)=0$ for any $k\geq 1$.
\end{theorem}

\begin{remark}
The mapping torus 
$M_\varphi$ for $\varphi\in\M_{g,1}$ as above 
admits no hyperbolic structures, 
so that 
$\mathrm{vol}(M_\varphi)=0$ holds. 
\end{remark}

\begin{proof}
At first, 
we prove the theorem for the genus one case. 
After that 
we give the outline of 
the proof in the higher genus case. 

Let 
$D_c$ be a Dehn twist 
along a non-separating simple closed curve $c$ 
on $\Sigma_{1,1}$. 
Taking a conjugation, 
we can assume that 
the curve $c$ is one of the standard generators 
of $\pi_1(\Sigma_{1,1})$. 
We then see that 
$\varphi={D_c}^q$ is represented by a matrix  
$\begin{pmatrix}
1&q\\
0&1
\end{pmatrix}
\in SL(2,\Z)$. 
Thus 
we can choose a deficiency one presentation 
$$
\left\langle 
x,y,t~|~txt^{-1}=x,~tyt^{-1}=x^qy \right\rangle
$$
of $\pi_1(M_\varphi)$. 
Applying the free differential calculus to 
the  relators 
$txt^{-1}x^{-1}$ and 
$tyt^{-1}(x^q y)^{-1}$, 
we obtain the Alexander matrix 
\[
A=\begin{pmatrix}
t-1 & 0 \\
-\partial(x^q)/\partial x & t-x^q 
\end{pmatrix}. 
\]
Here 
we remark that the generators 
$t$ and $x$ can be commuted 
by the relation $txt^{-1}=x$. 
Hence in this case, 
the $k$th Alexander matrix $A_k$ 
coincides with the original matrix $A$. 
In particular, 
 $t$ and $x$ always commute. 
As we saw in Section 5, 
the $L^2$-torsion invariant $\tau_k(M_\varphi)~(k\geq1)$ 
can be computed by using the usual determinant and the Mahler measure 
in such a situation. 
Since 
$$
\det A=(t-1)(t-x^q)
$$ 
is a generalized cyclotomic polynomial, 
we obtain 
$\log\tau_k(M_\varphi)=0$ 
as desired 
(see Corollary \ref{cor.6.2}). 

In the higher genus case, 
we can assume that the mapping class $\varphi$ is given by 
\begin{align*}
&\varphi_*(x_1)=x_1,~\varphi_*(x_2)=x_1^{q_1}x_2,
\ldots,\\
&\varphi_*(x_{2l-1})=x_{2l-1},~
\varphi_*(x_{2l})=x_{2l-1}^{q_l}x_{2l},\\
&\varphi_*(x_{2l+1})=x_{2l+1},\ldots,\varphi_*(x_{2g})=x_{2g}
\end{align*}
by taking a conjugation, 
where 
$q_1,\dots,q_l\in\Z$ 
and $1\leq l\leq g$. 
We then 
obtain the following presentation 
of $\pi_1(M_\varphi)$:
$$
\langle 
x_1,\dots,x_{2g},t~|~
tx_it^{-1}=\varphi_*(x_i),~ 
1\leq i\leq 2g
\rangle.
$$
Since 
the Alexander matrix  $A$ is the direct sum of 
the $2\times2$-matrix in the genus one case, 
so that 
we obtain $\log\tau_k(M_\varphi)=0$ 
by the similar arguments. 
\end{proof}

As another affirmative answer to Problem \ref{conj.3.7}, 
we can show the vanishing of $\log\tau_k$ 
for the following mapping classes 
(see \cite{KMT04-1}). 
That is, 
we consider the case 
where 
there exists an integer $n$ 
such that $M_{\varphi^n}$ is topologically the product 
of $\Sigma_{g,1}$ and $S^1$. 
Here 
its bundle structure is non-trivial 
in general. 
Namely 
the $n$th power $\varphi^n$ 
of a given monodromy $\varphi$ 
is not trivial. 
A typical example 
is the Dehn twist along the simple closed curve 
on $\Sigma_{g,1}$ 
parallel to the boundary. 
The difference between 
an isotopy fixing the boundary pointwisely 
and such one setwisely, 
it gives birth to the difference 
between a bundle structure and a topological type. 
We then obtain 

\begin{theorem}[\cite{KMT04-1}]
$\log\tau_k(M_{\varphi})=0$ for any $k\geq 1$. 
\end{theorem}

It is easy to see 
that such a 3-manifold does not admit a hyperbolic structure. 
Hence 
it has trivial simplicial volume. 

The above two examples are both non-hyperbolic cases, 
so that we conclude the present paper 
with the following problem. 

\begin{problem}
Show 
$$
\lim_{k\to\infty}\log\tau_k(M_\varphi)
=\log\tau(M_\varphi)
$$
for a pseudo-Anosov diffeomorphism $\varphi$. 
\end{problem}

\noindent
{\bf Acknowledgements.} 
The authors are grateful to the referee 
for his/her numerous and helpful comments 
which greatly improved this paper.


\end{document}